\documentstyle{article}
\newfont{\sm}{cmr10}                  \def\bij{b_{ij}}
\def\a{\alpha}                        \def\d{\delta}
\def\b{\beta}                         \def\D{\Delta}
                  \def\g{\gamma}
                       \def\gg{{\bf g}}
\def\L{\Lambda}

 \def\bg{{\bf g}}                     \def\R{{\cal R}}
\def\DY2{\widehat{DY(sl_2)}}         \def\Y-{\widehat{Y}^-}
\def\Yp{\widehat{Y}^+}             \def\ve{\varepsilon}
\def\Ypg{\widehat{Y}^+(\gg)}
\def\Ym{\widehat{Y}^-}
\def\Ymg{\widehat{Y}^-(\gg)}
\def\Ypmg{\widehat{Y}^\pm(\gg)}
  \def\semidir{\cdot}
\def\h{{\rm\hbar}}
\def\da{{\partial_\alpha}}
\def\nn{\nonumber}
\def\DYg{\widehat{DY({\bf g})}}
\def\fract#1#2{{\hbox{\footnotesize $#1$}\over\hbox{\footnotesize $#2$}}}

\def\q2a{q^{2(\a,\a)}}                \def\q-a{q^{-(\a,\a)}}
          
\newcommand{\bn}{\begin{equation}}     \newcommand{\ed}{\end{equation}}
\newtheorem{definition}{Definition}[section]
\newtheorem{proposition}{Proposition}[section]
\newtheorem{theorem}{Theorem}[section]

\newcommand{\hsp}{\mbox{$\hspace{.3in}$}}

\begin{document}
\begin{center}
\vspace{1cm}
{\LARGE\bf Central Extension of the Yangian Double}\\
\vspace{0.3cm}
{\large\bf S.M. Khoroshkin}\footnote{
e-mail:  sergejk@math.kth.se}
\\
\vspace{0.2cm}
 Institute of Theoretical and Experimental Physics,
117259 Moscow, Russia\\
\end{center}
\date{}
\vspace{0.2cm}
\begin{abstract}
Central extension $\DYg$ of the Double of the Yangian is defined for a
 simple Lie algebra ${\bf g}$ with complete proof for ${\bf g} =sl_2$.
 Basic representations and intertwining operators are constructed for
	$\DY2$.
\end{abstract}
\setcounter{footnote}{0}
\setcounter{equation}{0}

\section{Introduction}
The Yangian $Y({\bf g})$ was introduced by V. Drinfeld \cite{D1} as a Hopf
 algebra  quantizing the rational solution $r(u)=\frac{e_i\otimes e^i}{u}$
	of classical Yang-Baxter equation. As a Hopf algebra $Y(g)$ is a
	deformation  of the universal enveloping algebra $U(\bg[ u])$
	of polynomial currents to a simple Lie algebra $\bg$ with respect to
	cobracket defined by $r(u)$. Unfortunately, up to the moment the
 representation theory of the Yangian is not so rich in applications as
	it takes place, for instance, for quantum affine algebras \cite{JM}.
	 We could mention the following gaps:

	{\it (i)} The Yangian is not quasitriangular Hopf algebra, but
	pseudoquasitriangular Hopf algebra \cite{D1};

	{\it (ii)} There are no nontrivial examples of infinite-dimensional
	 representations of $Y(\bg)$.

	 In order to get quasitriangular Hopf algebra one should introduce the
		quantum double $DY(\bg )$ of the Yangian. Detailed analisys of
	$DY(\bg)$ was done in \cite{KT} together with explicit description of
	 the universal $R$-matrix for $DY(\bg )$  (complete for $DY(sl_2)$
 and partial in general case).

 The most important examples of infinite
dimensional representations of (quantum) affine algebras appear for nonzero
 value of central charge. Analogously, in the case of the Yangian we could
 expect the appearance of infinite dimensional representations only after
 central extension of $DY(\bg )$. This program is realized in this paper
 with complete proof for $DY(sl_2)$: we give a description of central
 extension $\DY2$ and construct its basic representations in bozonized
 form. In general case we present a description of $\DYg$ without complete
  proof. In the forthcoming paper \cite{KLP} we demonstrate, following
         general scheme of \cite{XXZ}, how our construction produce the
          formulas for correlation functions in rational models \cite{S}.

 The central extension of $DY(\bg)$ could be constructed in two ways.
 In Faddev-\-Reshetikhin-\-Takhtajan approach \cite{FRT}
 we could describe $DY(\bg)$
 by the set of equations
 \bn
 R_{V,W}^\pm(u-v)L_V^\pm(u)L_W^\pm(v)=L_W^\pm(v)L_V^\pm(u)
 R_{V,W}^\pm(u-v),
 \label{1.1}
 \ed
 \bn
 R_{V,W}^+(u-v)L_V^+(u)L_W^-(v)=L_W^-(v)L_V^+(u)
 R_{V,W}^+(u-v)
 \label{1.2}
 \ed
 for matrix valued generating functions $L_V^\pm(u)$,
  $L_W^\pm(v)$ of $DY(\bg)$, where $V(u), W(v)$ are finite
 dimensional representations of $DY(\bg)$, $R_{V,W}^\pm$ are
 images of the universal $R$-matrix $\R$ and $(\R^{-1})^{21}$.
  Then, following \cite{RS},
 we can make a shift of a spectral parameter in $R$-matrix in equation
  (\ref{1.2}) by a central element:
 \bn
 R_{V,W}^+(u-v+c\h)L_V^+(u)L_W^-(v)=L_W^-(v)L_V^+(u)
 R_{V,W}^+(u-v).
 \label{1.3}
 \ed
 For the construction of representations we should then extract Drinfeld
 generators of $\DYg$ from $L$-operators (\ref{1.3}) via their Gauss
 decomposition \cite{DF}.

 We prefer another way, originaly used by V.Drinfeld \cite{D2} for his
 ''current'' description of quantum affine algebras. The properties of
 comultiplication for $Y(\bg)$ show that one can extend the Yangian
 $Y(\bg)$ (or its dual with opposite comultiplication
 $Y^0(\bg)$) to a new Hopf algebra, adding the
 derivative $d$ of a spectral parameter $u$. Alternativly, one can extend
 $Y(\bg)$ by automorphisms of shifts of $u$. The double of this extension
 is exactly what we want to find. Central element $c$ is dual to
 derivative $d$.

 The plan of the paper is as follows. First we remind the description of
 the Yangian $Y(sl_2)$ and of its quantum double from \cite{KT}. Then we
 construct the central extension $\DY2$, describe the structure of the
 universal $R$-matrix for $\DY2$ and translate our description into $L$-
 operator language. In section 5 we construct basic representation of
 $\DY2$ and in the last section we describe the structure of $\DYg$ in
 general case (without a proof).
\setcounter{footnote}{0}
\setcounter{equation}{0}

\section{$Y(sl_2)$ and its Quantum Double}
The Yangian $Y(sl_2)$ can be described as a Hopf algebra generated by the
 elements $e_k, h_k, f_k$, $k\geq 0$ subjected to the relations
$$
[h_{k},h_{l}]=0 \ ,\hsp [e_{k}, f_{l}]= h_{k+l}\ ,$$
 $$[h_{0},e_{l}] = 2e_{l}\ , \hsp
[h_{0},f_{l}] = -2f_{l}\ ,$$
 $$
[h_{k+1},e_{l}] - [h_{k},e_{l+1}]=
\h\{ h_{k},e_{l}\}\ ,$$
 $$[h_{k+1},f_{l}] - [h_{k},f_{l+1}]=
-\h\{h_{k},f_{l}\}\ ,$$
 $$[e_{k+1},e_{l}] - [e_{k},e_{l+1}]=
\h\{ e_{k},e_{l}\}\ ,$$
\bn
[f_{k+1},f_{l}] - [f_{k},f_{l+1}]=
-h\{f_{k},f_{l}\}\ ,
\label{2.1}
\ed
 where $\h$ is a parameter of the deformation, $\{ a,b\} =ab+ba$.
 The comultiplication and the antipode are uniquely defined by the relations
 $$\D(e_0)= e_0\otimes 1+1\otimes e_0, \hsp
 \D(h_0)= h_0\otimes 1+1\otimes h_0, \hsp
 \D(f_0)= f_0\otimes 1+1\otimes f_0,$$
 $$\D(e_1)= e_1\otimes 1+1\otimes e_1 +\h h_0\otimes e_0,\hsp
 \D(f_1)= f_1\otimes 1+1\otimes f_1 +\h f_0\otimes h_0,$$
 \bn
 \D(h_1)= h_1\otimes 1+1\otimes h_1+\h h_0\otimes h_0-2\h f_0
\otimes e_0.
 \label{2.2}
 \ed
 In terms of generating functions
 $$ e^{+}(u):=\sum_{k\geq
0}e_{k}u^{-k-1} , \hsp f^{+}(u):=\sum_{k\geq
0}f_{k}u^{-k-1}, \hsp h^{+}(u):=1+\h\sum_{k\geq
0}h_{k}u^{-k-1}$$
the relations (\ref{2.1}) look as follows
$$
[h^{+}(u),h^{+}(v)]=0\ ,$$
$$
[e^{+}(u),f^{+}(v)]=-\frac{1}{\h}\frac{h^{+}(u)-h^{+}(v)}{u-v}\ ,$$
$$
[h^{+}(u),e^{+}(v)]=
-\h\frac{\{h^{+}(u),(e^{+}(u)-e^{+}(v))\}}{u-v}\ ,$$
$$
[h^{+}(u),f^{+}(v)]=
\h\frac{\{h^{+}(u),(f^{+}(u)-f^{+}(v))\}}{u-v}\ ,$$
$$
[e^{+}(u),e^{+}(v)]=
-\h\frac{(e^{+}(u)-e^{+}(v))^2}{u-v}\ ,$$
\bn
[f^{+}(u),f^{+}(v)]=
\h\frac{(f^{+}(u)-f^{+}(v))^2}{u-v}\ .
\label{2.3}
\ed
The comultiplication is given by Molev's formulas \cite{M}, see also
\cite{KT}:
\bn
\D(e^{+}(u))=e^{+}(u)\otimes 1 +\sum_{k=0}^\infty (-1)^k\h^{2k}
\big(f^{+}(u+\h)\big)^kh^{+}(u)\otimes\big(e^{+}(u)\big)^{k+1},
\label{2.4}
\ed
\bn
\D(f^{+}(u))=1\otimes f^{+}(u)+\sum_{k=0}^\infty (-1)^k\h^{2k}
\big(f^{+}(u)\big)^{k+1}\otimes h^{+}(u)\big(e^{+}(u+\h)\big)^k
\label{2.5}
\ed
\bn
\D(h^{+}(u))= \sum_{k=0}^\infty (-1)^k(k+1)\h^{2k}
\big(f^{+}(u+\h)\big)^kh^{+}(u)\otimes h^{+}(u)
\big(e^{+}(u+\h)\big)^{k}
\label{2.6}
\ed

Let now $C$ be an algebra generated by the elements $e_{k}$, $f_{k}$,
$h_{k}$, $(k\in {\bf Z})$, with relations
(\ref{2.1}). Algebra $C$ admits {\bf Z}-filtration
\bn
\ldots\subset C_{-n}\subset \ldots\subset C_{-1} \subset C_{0}
\subset C_{1}\ldots\subset C_{n}\ldots \subset C
\label{2.7}\label{completion}
\ed
defined by the conditions $\mbox{deg}\, e_{k}=\mbox{deg}\, f_{k}=
\mbox{deg}\, h_{k}=k$; $\deg x\in C_m\leq m$. Let $\bar{C}$ be the
corresponding formal completion of $C$. It is proved in \cite{KT}
that $DY(sl_2)$ is isomorphic to $\bar{C}$ as an algebra.
 In terms of generating functions
 $$ e^{\pm}(u):=\pm\sum_{k\geq
0\atop k<0}e_{k}u^{-k-1} , \hsp f^{\pm}(u):=\pm\sum_{k\geq
0\atop k<0}f_{k}u^{-k-1}, \hsp h^{\pm}(u):=1\pm\h\sum_{k\geq
0\atop k<0}h_{k}u^{-k-1},$$
 $$e(u)=e^+(u)-e^-(u),\hsp f(u)=f^+(u)-f^-(u)$$
 the defining relations for $DY(sl_2)$ look as follows:
\begin{eqnarray}
h^a(u)h^b(v)&=&h^b(v)h^a(u), \hsp a,b=\pm, \nn\\
e(u)e(v)&=&{u-v+\h\over u-v-\h}\ e(v)e(u) \nn\\
f(u)f(v)&=&{u-v-\h\over u-v+\h}\ f(v)f(u) \nn\\
h^\pm(u)e(v)&=&{u-v+\h\over u-v-\h}\ e(v)h^\pm(u) \nn\\
h^\pm(u)f(v)&=&{u-v-\h\over u-v+\h}\ f(v)h^\pm(u) \nn\\
{[}e(u),f(v){]}&=&
{1\over \h}\delta(u-v)\left(h^+(u)- h^-(v)\right)
\label{comDY}
\end{eqnarray}
Here
$$\delta(u-v)=\sum_{n+m=-1}u^nv^m$$
satisfies the property
$$\delta(u-v)f(u)=\delta(u-v)f(v).$$
The comultiplication is given by the same formulas
$$
\D(e^{\pm}(u))=e^{\pm}(u)\otimes 1 +\sum_{k=0}^\infty (-1)^k\h^{2k}
\big(f^{\pm}(u+\h)\big)^kh^{\pm}(u)\otimes\big(e^{\pm}(u)\big)^{k+1},$$
$$\D(f^{\pm}(u))=1\otimes f^{\pm}(u)+\sum_{k=0}^\infty (-1)^k\h^{2k}
\big(f^{\pm}(u)\big)^{k+1}\otimes h^{\pm}(u)\big(e^{\pm}(u+\h)\big)^k
$$
\bn
\D(h^{\pm}(u))= \sum_{k=0}^\infty (-1)^k(k+1)\h^{2k}
\big(f^{\pm}(u+\h)\big)^kh^{\pm}(u)\otimes h^{\pm}(u)
\big(e^{\pm}(u+\h)\big)^{k}.
\label{Molev}
\ed
Subalgebra $Y^+=Y(sl_2)\subset DY(sl_2)$ is generated by the
components of $e^+(u)$, $f^+(u)$, $h^+(u)$ and its dual with
opposite comultiplication $Y^-=(Y(sl_2))^0$ is a formal
completion (see (\ref{completion})) of subalgebra, generated by
 the components of $e^-(u)$, $f^-(u)$, $h^-(u)$.

 Let us describe the Hopf pairing of $Y^+$ and $Y^-$ \cite{KT}.
 Note that by a Hopf pairing $<,>: A\otimes B \rightarrow {\bf
 C}$ we mean bilinear map satisfying the conditions
 \bn
 <a, b_1b_2>=<\D(a), b_1\otimes b_2>, \hsp <a_1a_2,b>=<a_2\otimes
 a_1, \D(b)> \label{Hopf}
 \ed
 The last condition is unusual but convenient in a work with
  quantum double.

 Let $E^\pm$, $H^\pm$, $F^\pm$ be subalgebras  (or
 their completions in $Y^-$ case) generated by the components of
  $e^\pm(u)$, $h^\pm(u)$, $f^\pm(u)$. Subalgebras $E^\pm$ and
 $F^\pm$ do not contain the unit.

The first property of the Hopf pairing $Y^+\otimes Y^-
\rightarrow {\bf C}$ is that it
 preserves the  decompositions
 $$ Y^+=E^+H^+F^+,\hsp Y^-=F^-H^-E^-.$$
 It means that
 \bn
 <e^+h^+f^+,f^-h^-e^->= <e^+,f^-><h^+,h^-><f^+,e^->
 \label{factorization}
 \ed
 for any elements $e^\pm\in E^\pm$, $h^\pm\in H^\pm$,
 $f^\pm\in F^\pm$.
 This property defines the pairing uniquely together with the
 relations
 $$<e^+(u),f^-(v)>= \frac{1}{\h (u-v)},\hsp
 <f^+(u),e^-(v)>= \frac{1}{\h (u-v)},$$
 $$<h^+(u),h^-(v)>=\frac{u-v+\h}{u-v-\h}.$$
 The full information of the pairing is encoded in the universal
  $R$-matrix for $DY(sl_2)$ which has the following form
\cite{KT}:  \bn \R=R_+R_0R_- \ , \label{UR29} \ed where \bn
R_+=\prod_{k\geq 0}^{\rightarrow}\exp(-\h e_k\otimes f_{-k-1})\ ,
\hsp R_-=\prod_{k\geq 0}^{\leftarrow}\exp(-\h f_k\otimes
e_{-k-1})\ , \label{UR30} \ed \bn R_0=\prod_{n\geq 0}\exp\h{\rm
Res}_{u=v}(\frac{\rm d}{{\rm d}u}k^+(u))\otimes k^{-}(v+2n\h
+\h), \label{UR31} \ed Here $k^{\pm}(u)=\frac{1}{\h}\ln
h^{\pm}(u)$.

\setcounter{footnote}{0}
\setcounter{equation}{0}
\section{Central extension of $DY(sl_2)$}
 Let $d=\frac{{\rm d}}{{\rm d}u}$ be the operator of derivation
of a spectral parameter: $d\, g(u) =\frac{{\rm d}}{{\rm d}u}
g(u)$ and $T_x = \exp(xd)$ be the shift operator: $T_xg(u) =
g(u+x)$.

Let us define semidirect products $Y^\pm\semidir{\bf C}[[d]]$,
 $\mbox{deg}\, d = -1$,
in a
natural way:
\bn
[d,e^\pm(u)]=\frac{{\rm d}}{{\rm d}u}e^\pm(u), \hsp
[d,h^\pm(u)]=\frac{{\rm d}}{{\rm d}u}h^\pm(u), \hsp
[d,f^\pm(u)]=\frac{{\rm d}}{{\rm d}u}f^\pm(u).
\label{3.0}
\ed
\begin{proposition} \label{proposition 3.1}
Semidirect products $Y^\pm\semidir{\bf C}[[d]]$ are Hopf algebras
 if we put
$$\D(d)=d\otimes 1 +1\otimes d$$
\end{proposition}
The proof follows by induction from the observation that
$$\D(a_i)=a_i\otimes 1+1\otimes a_i +\;{\mbox{ terms of degree
lower then }}i$$
 for $a=e,h,f$. Another argument is that the coproduct of
$a^\pm(u)$ can be expressed again in terms of $a^\pm(u)$.

Denote by $\Ym$ the Hopf algebra $Y^-\semidir{\bf C}[[d]]$.
 Let $\Yp$ be the following Hopf algebra: $\Yp$ is a tensor
product of $Y^+$ and of polynomial ring of central element $c$,
 $\mbox{deg}\, c = 0$, as
an algebra:  $$\Yp=Y^+\otimes{\bf C}[c],$$
$$\D_{\Yp}(c)=c\otimes 1+1\otimes c$$ and
\bn
\D_{\Yp}a^+(u)=\left({\rm Id}\otimes T_{\h c\otimes 1}^{-1}
\right) \D_{Y^+}a^+(u) \label{3.1} \ed for $a= e,h,f$. For
instance, $$\D_{\Yp}e^+(u)= e^{+}(u)\otimes 1 +\sum_{k=0}^\infty
(-1)^k\h^{2k}
\big(f^{+}(u+\h)\big)^kh^{+}(u)\otimes\big(e^{+}(u-\h
c_1)\big)^{k+1},$$
   where $c_1=c\otimes 1$.
\begin{proposition} \label{proposition 3.2}
There exists unique extension of the Hopf pairing from
 $Y^p\otimes Y^m$ to $\Yp\otimes \Ym$ satisfying the conditions:

{\em (i)} $\hsp <c,d>=\frac{1}{\h};$

{\em (ii)} The Hopf pairing preserves the decompositions
$$\Yp =Y^+{\bf C}[c], \hsp   \Ym =Y^-{\bf C}[[d]]$$
\end{proposition}
{\bf Proof.} Let $a^+(u)\in Y^+$, $b^-(u)\in Y^-$,
$<a^+(u),b^-(v)>= f(u-v)$ and $\g$ is a number, $\g\in{\bf C}$.
Then, due to (i), (ii) and (\ref{3.1})
$$ <a^+(u), e^{\g d}b^-(v)e^{-\g d}>=
 <a^+(u), e^{\g d}b^-(v)>=$$
$$ < \D_{\Yp}a^+(u), e^{\g d}\otimes b^-(v)>  =f(u-x-\g)=
<a^+(u), b^-(v+\g)>.$$
which proves the compatibility of the extended pairing with the
relations (\ref{3.0}).
\begin{definition}
Central extension $\DY2$ of $DY(sl_2)$ is quantum double
$D
(\Yp)$ of the Hopf algebra $\Yp =\Yp(sl_2)$.
\label{dfn3.1}
\end{definition}
Equivalently, $\DY2$ is the double $D(\Ym)$ with opposite comultiplication.

The following theorem describes $\DY2$ explicitely as a Hopf algebra.
\begin{theorem}
$\DY2$ is isomorphic to a formal completion (see (\ref{2.7}))
   of
 the algebra with generators $e_k, f_k, h_k$, $k\in {\bf Z}$, $d$ and central
 element $c$ with the relations written in terms of generating functions:
$$[d,e(u)]=\frac{\mbox{d}}{\mbox{d}u}e(u), \hsp
 [d,f(u)]=\frac{\mbox{d}}{\mbox{d}u}f(u), \hsp
 [d,h^\pm(u)]=\frac{\mbox{d}}{\mbox{d}u}h^\pm(u),$$
\begin{eqnarray}
e(u)e(v)&=&{u-v+\h\over u-v-\h}\ e(v)e(u) \nn\\
f(u)f(v)&=&{u-v-\h\over u-v+\h}\ f(v)f(u) \nn\\
h^\pm(u)e(v)&=&{u-v+\h\over u-v-\h}\ e(v)h^\pm(u) \nn\\
h^+(u)f(v)&=&{u-v-\h-\h c\over u-v+\h-\h c}\ f(v)h^+(u) \nn\\
h^-(u)f(v)&=&{u-v-\h\over u-v+\h}\ f(v)h^-(u) \nn\\
h^+(u)h^-(v)&=&{u-v+\h\over u-v-\h}\cdot{u-v-\h-\h c\over u-v+\h-\h c}\
h^-(v)h^+(u)
\nn\\
{[}e(u),f(v){]}&=&
{1\over \h}\left(\delta(u-(v+\h c))h^+(u)-\delta(u-v)h^-(v) \right)
\label{3.5}
\end{eqnarray}
The comultiplication is given by the relations
$$
\D(e^{\varepsilon}(u))=e^{\ve}(u)\otimes 1
+\sum_{k=0}^\infty (-1)^k\h^{2k}
\big(f^{\ve}(u+\h-\d_{\varepsilon ,+}\h c_1)\big)^k
h^{\ve}(u)\otimes\big(e^{\ve}(u-\d_{\varepsilon ,+}\h c_1)
\big)^{k+1},$$
 $$\D(h^{\varepsilon}(u))=
\sum_{k=0}^\infty (-1)^k(k+1)\h^{2k}
\big(f^{\ve}(u+\h-\d_{\varepsilon ,+}\h c_1)\big)^kh^{\ve}(u)\otimes
h^{\ve}(u-\d_{\varepsilon ,+}\h c_1) \big(e^{\ve}(u+\h-
\d_{\varepsilon ,+}\h c_1)\big)^{k},$$
\bn
\D(f^{\varepsilon}(u))=1\otimes f^{\ve}(u)
+\sum_{k=0}^\infty (-1)^k\h^{2k}
\big(f^{\ve}(u+\d_{\varepsilon ,+}\h c_2)\big)^{k+1}\otimes
h^{\ve}(u)\big(e^{\ve}(u+\h)\big)^k.
 \label{3.6}
 \ed
 where  $\varepsilon =\pm,\; \d_{+ ,+}=1$ and $\d_{- ,+}=0$.
 \label{th3.1}
 \end{theorem}
 The proof of Theorem
 \ref{th3.1} reduces to explicit calculation of commutation relations
 in quantum double of $\Yp$. In abstract Sweedler notation for a
 double of a Hopf algebra $A$ these relations have the following form:
 \bn
 a\cdot b=<a^{(1)},b^{(1)}>\,<S^{-1}(a^{(3)}),b^{(3)}>\,
 b^{(2)}\cdot a^{(2)}
 \label{3.7}
 \ed
 where
 $a\in A,\;b\in A^0,\; \D^2(x)=(\D\otimes Id)\D(x)=x^{(1)}\otimes
 x^{(2)}\otimes x^{(3)}$, $S$ is  antipode in $A$.

 For the calculation of (\ref{3.7}) we need the following partial
 information about $\D^2$ and $S^{-1}$ in $\widehat{Y}^\pm$ which
  one can deduce from (\ref{Molev}) or directly by induction:
 $$ S^{-1}e^+(u)=-e(u-\h c)h^{-1}(u-\h c) \;\;\;{\rm mod}\; E^+\Yp F^+,$$
 $$ S^{-1}f^+(u)=-h^{-1}(u-\h c)f(u-\h c) \;\;\;{\rm mod}\; E^+\Yp F^+,$$
 \bn
  S^{-1}h^+(u)=h^{-1}(u-\h c) \;\;\;{\rm mod}\; E^+\Yp F^+
  \label{3.8}
  \ed
  and

  $$\D^2e^\pm(u)=K^\pm\left(e^\pm_1(u) +
  h^\pm_1(u) e^\pm_2(u) +
  h^\pm_1(u) h^\pm_2(u) e^\pm_3(u)-\right.$$
 $$\left. - \h^2h^\pm_1(u)f^\pm_1(u-\h)(e^\pm_2(u))^2-
 2\h^2h^\pm_1(u)f^\pm_1(u-\h) e^\pm_2(u-\h)h^\pm_2(u)
 e^\pm_3(u)\right)
 \;\;\;{\rm mod}\;X^\pm$$

 $$\D^2f^\pm(u)=K^\pm\left(f^\pm_3(u) +
  f^\pm_2(u)h^\pm_3(u) +
 f^\pm_1(u)h^\pm_2(u) h^\pm_3(u) -\right.$$
 $$\left. - \h^2(f^\pm_2(u))^2e^\pm_3(u-\h)h^\pm_3(u)-
	 2\h^2f^\pm_1(u)h^\pm_2(u)f^\pm_2(u-\h) e^\pm_3(u-\h)h^\pm_3(u)
	 \right)
	 \;\;\;{\rm mod}\;X^\pm$$

 $$ \D^2h^\pm(u)=K^\pm\left( h_1^\pm(u)h_2^\pm(u)h^\pm_3(u)-
	 2\h^2h^\pm_1(u)h^\pm_2(u)f^\pm_2(u-\h) e^\pm_3(u-\h)h^\pm_3(u)-
	 \right.$$
	 \bn
	 \left.
	 -2\h^2h^\pm_1(u)f^\pm_1(u-\h) e^\pm_2(u-\h)h^\pm_2(u)
	 h^\pm_3(u)\right)
	 \;\;\;{\rm mod}\;X^\pm
	 \label{3.9}
	 \ed
	 where
	 $$K^+=Id \otimes T_{\h c_1}^{-1}\otimes T_{\h(c_1+c_2)}^{-1},
	 \hsp K^-=Id\otimes Id\otimes Id,$$
	 $$X^\pm =(H^\pm F^\pm +H^\pm)C^\pm\otimes\widehat{Y}^\pm
	 \otimes E^{\pm 2}C^\pm +F^{\pm 2}C^\pm\otimes\widehat{Y}^\pm
	 \otimes (E^\pm H^\pm  +H^\pm)C^\pm,$$
	 $C^+ ={\bf C}[c], C^- ={\bf C}[[d]]$, and, as usually, $a_1(u)$
	 means $a(u) \otimes 1\otimes 1$, $a_2(u) = 1\otimes a(u)\otimes 1$,
	  $a_3(u)= 1\otimes 1\otimes a(u)$.

		The substitution of (\ref{3.9}) into (\ref{3.8}) gives the
 following relations between generators of $\Yp$ and $\Ym$:
$$
[e^{+}(u),e^{-}(v)]=
-\h\frac{(e^{+}(u)-e^{-}(v))^2}{u-v}\ ,$$
$$
[e^{+}(u),f^{-}(v)]=-\frac{1}{\h}\frac{h^{+}(u)-h^{-}(v)}{u-v}\ ,$$
$$
[e^{+}(v),h^{+}(v)]=
 \h\frac{\{h^{-}(v),(e^{+}(u)-e^{-}(v))\}}{u-v}\ ,$$

$$
[h^{+}(u),e^{-}(v)]=
-\h\frac{\{h^{+}(u),(e^{+}(u)-e^{-}(v))\}}{u-v}\ ,$$
$$
[h^{+}(u),f^{-}(v)]=
\h\frac{\{h^{+}(u),(f^{+}(u)-f^{-}(v))\}}{u-v}\ ,$$
$$h^+(u)h^-(v)={u-v+\h\over u-v-\h}\cdot{u-v-\h-\h c\over u-v+\h-\h c}\
h^-(v)h^+(u)$$
$$
[f^{-}(u),e^{+}(v)]=-\frac{1}{\h}\left(\frac{h^{+}(u)}{u-v}-
\frac{h^{-}(v)}{u-v-\h c}\right) ,$$

$$
[f^{+}(u),h^{-}(v)]=
\h\frac{\{h^{-}(v),(f^{-}(v)-f^{+}(u))\}}{u-v-\h c}\ ,$$
$$
[f^{+}(u),f^{-}(v)]=
\h\frac{(f^{+}(u)-f^{-}(v))^2}{u-v}$$
$$
e^+(u)e^{\g d}=e^{\g d}e^+(u-\g),\hsp
h^+(u)e^{\g d}=e^{\g d}h^+(u-\g),$$
\bn
f^+(u)e^{\g d}=e^{\g d}f^+(u-\g).
\label{3.10}
\ed
After a change of variables
$$f^+_{\rm new}(u)=f^+_{\rm old}(u+\h c),$$
$$f(u)=f^+_{\rm new}(u)-f^-(u),\hsp e(u)=e^+(u)-e^-(u)$$
we get the relations (\ref{3.5}) and comultiplication rules
 (\ref{3.6}). The theorem is proved.

 As a consequence of Proposition \ref{proposition 3.2} and of the
 description of the universal $R$-matrix for $DY(sl_2)$ we have the
	explicit formula for the universal $R$-matrix for $\DY2$.

	\begin{theorem}
  \bn \R=R_+R_0R_-\exp(\h c\otimes d) \ ,
	\label{5.29} \ed
	where
$$
R_+=\prod_{k\geq 0}^{\rightarrow}\exp(-\h e_k\otimes f_{-k-1})\ ,
\hsp R_-=\prod_{k\geq 0}^{\leftarrow}\exp(-\h g_k\otimes
e_{-k-1})\ ,$$
\bn
g_k=\sum_{m=0}^{k}\frac{k!}{m!(k-m)!}f_{k-m}(\h c)^m,
\label{5.30}
\ed
\bn
R_0=\prod_{n\geq 0}\exp\h{\rm
Res}_{u=v}(\frac{\rm d}{{\rm d}u}k^+(u))\otimes k^{-}(v+2n\h
+\h),
\label{5.31}
\ed
Here $k^{\pm}(u)=\frac{1}{\h}\ln
h^{\pm}(u)$.
\label{th3.2}
\end{theorem}
The universal $R$-matrix for $\DY2$ can be rewritten also in slightly
 more symmetric form:
$$ \R =R'_+R'_0R'_-$$
where
$$R'_+=R_+,\hsp R'_0=R_0 \exp(\h c\otimes d),\hsp
R'_-=\prod_{k\geq 0}^{\leftarrow}\exp(-\h f_k\otimes
e_{-k-1})\ .$$

\setcounter{footnote}{0}
\setcounter{equation}{0}
\section{$L$-operator presentation of $\DY2$}
Let $\rho(x) = T_{-x}\rho(0)$ be the action of $\DY2$ in two-dimensional
 evaluation representation $W_x$ of $DY(sl_2)$. In a basis $w_\pm$ this
	action looks as follows:
 \bn
 e_kw_+=f_kw_-=0,\; e_kw_-=x^kw_+,\;f_kw_+=x^kw_-,\; h_kw_+=x^kw_+,\;
 h_kw_-=-x^kw_-
 \label{4.1}
 \ed
 Let
 $$L^+(x)=(\rho(x)\otimes Id)\exp(\h d\otimes c)(\R^{21})^{-1},$$
 $$L^-(x)=(\rho(x)\otimes Id)\R\exp(-\h c\otimes d),$$
 and
 $$ R^+(x-y)=(\rho(x)\otimes \rho(y))\exp(\h d\otimes c)(\R^{21})^{-1},$$
 $$R^-(x-y)=(\rho(x)\otimes \rho(y))\R\exp(-\h c\otimes d).$$
 We have \cite{KT}
 \bn
 R^\pm(u)=\rho^\pm(u)\cdot\left(
 \begin{array}{cccc}
 1&0&0&0\\
 0&\frac{u}{u+\h}&\frac{\h}{u+\h}&0\\
 0&\frac{\h}{u+\h}&\frac{u}{u+\h}&0\\
 0&0&0&1
 \end{array}
 \right)
 \label{4.2}
 \ed
 where
 $$\rho^\pm(u)=\left(
        \frac{\Gamma(\mp\frac{u}{2\h})\Gamma(1\mp\frac{u}{2\h})}
 {\Gamma^2(\frac{1}{2}\mp\frac{u}{2\h})}\right)^{\pm 1}$$
 The Yang-Baxter equation on $\R$ implies the following relations on $L^+$
	and $L^-$ (see \cite{FR} for detailes in $U_q(\widehat{\gg})$ case):
 \bn
 R_{12}^\pm(x-y)L_1^{\pm}(x)L_2^\pm(y)=L_2^\pm(y)L_1^{\pm}(x)
 R_{12}^\pm(x-y)
 \label{4.4}
 \ed
 \bn
 R_{12}^+(x-y-\h c)L_1^{-}(x)L_2^+(y)=L_2^+(y)L_1^{-}(x)
 R_{12}^+(x-y)
 \label{4.5}
 \ed
where $L_1=L\otimes Id$, $L_2=Id\otimes L$. The properties of
 comultiplication for $\R$:
$$(\D\otimes Id)\R =\R^{13}\R^{23}, \hsp
	(Id\otimes\D)\R =\R^{13}\R^{12}$$
	imply the comultiplication rules for $L^\pm$:
        $${\D}' L^+(x)=L^+(x-\h c_2)\otimes L^+(x),\hsp{\rm or}\;\;\;\;
	 \D l_{ij}^+(u)=\sum_kl^+_{kj}(u)\otimes l_{ik}^+(u-\h c_1),$$
        $${\D}' L^-(x)=L^-(x)\otimes L^-(x),\hsp{\rm or}\;\;\;\;
	 \D l_{ij}^-(u)=\sum_kl_{kj}^-(u)\otimes l_{ik}^-(u).$$
	 The explicit formula for the universal $R$-matrix and (\ref{4.1})
	 express Gauss factors of the $L$-operators in terms of Drinfeld
	 generators:
	 $$L^+(x)=
	 \left(\begin{array}{cc}1&\h f^+(x-\h c)\\0&1\end{array}\right)
	 \left(\begin{array}{cc}k_1^+(x)&0\\0&k_2^+(x)\end{array}\right)
	 \left(\begin{array}{cc}1&0\\\h e^+(x)&1\end{array}\right),$$

	 $$L^-(x)=
	 \left(\begin{array}{cc}1&\h f^-(x)\\0&1\end{array}\right)
	 \left(\begin{array}{cc}k_1^-(x)&0\\0&k_2^-(x)\end{array}\right)
	 \left(\begin{array}{cc}1&0\\\h e^-(x)&1\end{array}\right),$$
	 with
	 $h^\pm(x)= k_2^\pm(x)^{-1}k_1(x),$ $\;\;k_1^\pm(x)k_2^\pm(x-\h)
	 =1$ analogous to Ding-Frenkel formulas for $U_q(\widehat{gl}_n)$
	  \cite{DF}.

  Note that Frenkel and Reshetikhin use more symmetric form of the
  equation (\ref{4.5}), using the shifts of spectral parameter in both
  sides of equation. One can get analogous form via twisting of
  comultiplication in $\DY2$ by elements $\exp(\a\h c\otimes d)$ and
   $\exp(\b\h d\otimes c)$, $\a ,\b\in{\bf C}$. For quantum affine
   algebras such renormalizations are initiated by the condition of
   consistency with Cartan involution in quantum affine algebra.
   There is no analogous motivation in our case so we do not use
   these twists.

\setcounter{footnote}{0}
\setcounter{equation}{0}
\section{Basic representations of $\DY2$}
Let ${\cal H}$ be Heisenberg algebra generated by free bosons with
 zero mode $a_{\pm n}, n\geq 1, a_0= \a$ and $p= \partial_\a$ with
	commutation relations
	\bn
	[a_n,a_m]=n\delta_{n+m,0},\hsp [p,a_0]=2.
	\label{5.1}
	\ed
	In the following we use the generating functions
	\bn
	a_+(z)=\sum_{n\geq 1}\frac{a_n}{n}z^{-n}-p\log z,\hsp
	a_-(z)=\sum_{n\geq 1}\frac{a_{-n}}{n}z^{n}+\frac{a_0}{2},
	\label{5.2}
	\ed
	\bn
	a(z)=a_+(z)-a_-(z),\hsp \phi_\pm(z)=\exp a_\pm(z).
	\label{5.2a}
	\ed
	They satisfy the relation
	$$[a_+(z),a_-(u)]=-\log(z-w)$$

	Let $\Lambda_i, i=0,1$ be formal power extensions of the Fock spaces:
	$$\Lambda_i={\bf C}[[a_{-1},\ldots ,a_{-n},\ldots ]]\otimes
	\left(\oplus_{n\in{\bf Z}+\frac{i}{2}}{\bf C}e^{n\a}\right).$$
	The following relations define an action of $\DY2$ on $\Lambda_i$
	with central charge $c=1$. We call them basic representations of
	$\DY2$:
\begin{eqnarray}
e(u)&=&
\exp\left(\sum_{n=1}^{\infty} {a_{-n}\over n}\left[(u-\h)^n+u^n\right]\right)
\exp\left(-\sum_{n=1}^{\infty} {a_{n}\over n}u^{-n} \right)
e^{\alpha}u^{\da} \nn\\
f(u)&=&
\exp\left(-\sum_{n=1}^{\infty} {a_{-n}\over n}\left[(u+\h)^n+u^n\right]\right)
\exp\left(\sum_{n=1}^{\infty} {a_{n}\over n}u^{-n} \right)
e^{-\alpha}u^{-\da} \nn\\
h^-(u)&=&
\exp\left(\sum_{n=1}^{\infty} {a_{-n}\over n}\left[(u-\h)^n-(u+\h)^n\right]
\right)\nn\\
h^+(u)&=&
\exp\left(\sum_{n=1}^{\infty} {a_{n}\over n}\left[(u-\h)^{-n}-u^{-n}\right]
\right)
\left({u\over u-\h}\right)^\da
                                \label{5.3}
\end{eqnarray}
The action of $e^{\g d}$, $\g\in {\bf C}$ is defined by the prescriptions
\bn
e^{\g d}\cdot (1\otimes 1)=1\otimes 1
\label{5.4}
\ed
and
$$e^{\g d}a_{-n}e^{-\g d}=\sum_{k\geq 0}\frac{(n+k-1)!}{(n-1)!k!}
a_{-(n+k)}\g^k, \hsp n\geq 1,$$
$$e^{\g d}a_{n}e^{-\g d}=\sum_{0\leq k< n}(-1)^{k}\frac{(n)!}{(n-k)!k!}
a_{n-k}\g^k+(-1)^n\g^np, \hsp n\geq 1,$$
\bn
e^{\g d}a_{0}e^{-\g d}=a_0+2\left(\sum_{n\geq1}\frac{a_{-n}}{n}\g^n
\right),\hsp e^{\g d}pe^{-\g d}= p
\label{5.6}
\ed
Note that the relations (\ref{5.6}) define an automorphism
 of Heisenberg algebra $T_\g a(z) =  a(z+\g)$.

 In terms of generating functions $\phi_\pm(z)$ the action of $\DY2$
 in basic representations has the following compact form:
 $$e(u)=\phi_-(u-\h)\phi_-(u)\phi_+^{-1}(u),$$
 $$f(u)=\phi_-^{-1}(u+\h)\phi_-^{-1}(u)\phi_+(u),$$
 $$h^+(u)=\phi_+(u-\h)\phi_+^{-1}(u),$$
 $$h^-(u)=\phi_-(u-\h)\phi_-^{-1}(u+\h),$$
 \bn
 e^{\g d}\phi_\pm(u)=\phi_\pm(u+\g)e^{\g d}
 \label{5.7}
 \ed

The relations (\ref{5.3}), (\ref{4.1}), (\ref{3.6}) give
 possibility to write down  bozonized expressions for intertwining
 operators of type one and two
 $\Phi(z): \L_i \rightarrow \L_{1-i}\otimes W_z$,
 $\Psi(z): \L_i \rightarrow W_z\otimes \L_{1-i}$ analogous to the case of
 $U_q(\widehat{sl}_2)$ \cite{XXZ}. More concretely, let
$$\Phi(z): \L_i\to \L_{1-i}\otimes W_z,\quad
\Psi(z): \L_i\to W_z\otimes \L_{1-i},\quad x\in\DY2$$
 satisfy the relations
$$\Phi(z) x = \Delta(x) \Phi(z),\quad
\Psi(z) x = \Delta(x) \Psi(z).$$
The components of intertwining operators are defined as follows
$$\Phi(z) v  =  \Phi_+(z)v\otimes v_+ + \Phi_-(z)v\otimes v_-,\quad
\Psi(z) v  =  v_+\otimes\Psi_+(z)v + v_-\otimes\Psi_-(z)v,$$
where $v\in \L_0$ or in $\L_1$.
Then, for instance, $\Psi_-(z)$ is the solution of the system of equations
$$\Phi_-(z)h^+(u)={u-z-2\h\over u-z-\h}h^+(u)\Phi_-(z) ,$$
\bn\Phi_-(z)h^-(u)={u-z-\h\over u-z}h^-(u)\Phi_-(z) ,\label{phiequations}\ed
$$\Phi_-(z)e(u)=e(u)\Phi_-(z) ,$$
and
$$\Phi_+(z)=\Phi_-(z)f_0-f_0\Phi_-(z).$$
Formal solution of (\ref{phiequations}) is
\bn
\Phi_-(z)=\phi_-(z+\h)(-1)^{p\over 2}\prod_{k=0}^{\infty}
\frac{\phi_+(z-\h-2k\h)}{\phi_+(z-2k\h)}
\label{phi}
\ed
modulo normalization factors depending on $i$.\\
S.Pakuliak suggested the following regularization of the formal solution
(\ref{phi}):
\begin{eqnarray}
\Phi^{}_-(z)&=&\exp\left(\sum_{n=1}^{\infty}{a_{-n}\over n}(z+\h)^n\right)
e^{\alpha/2}  (2\h)^{\partial_\a/2}
\left( {{\Gamma(\fract{1}{2}-\fract{z}{2\h})}\over{\Gamma(-\fract{z}{2\h})}}
  \right)^{\partial_\a}
\nn\\
&\times&\prod_{k=0}^{N}
\exp\left(-\sum_{n=1}^{\infty}  {a_n\over n}
\left[(z-2k\h)^{-n}-(z-\h-2k\h)^{-n}\right]    \right)
\label{Phi-}
\end{eqnarray}

 Let us note that
 our form of presentation (\ref{5.7})  for basic
 representations and (\ref{phi}) for corresponding intertwining operators
 $\Phi(z)$ and $\Psi(z)$ looks to be quite general. In particular,
 one can see that the Frenkel-Jing
 \cite{FJ} formulas for basic representations of $U_q(\widehat{sl}_2)$
 and corresponding expressions for intertwining operators $\Phi(z)$
	and $\Psi(z)$ \cite{JM} could be rewritten in a form analogous to
	(\ref{5.7}) and (\ref{phi}) using usual boson field $a(z)$ and
	 multiplicative shifts of spectral parameter.

 It will be interesting to equip basic representations of $\DY2$
 with invariant bilinear form and with a structure of topological
        representation.

        An application of intertwining operators for $\DY2$ to calculation
         of correlation functions in $su(2)$-invariant Thirring model is
         given in forthcoming paper \cite{KLP}.

\setcounter{footnote}{0}
\setcounter{equation}{0}
\section{The general case}
The arguments of section 3 prove that there exists central extension
 $\DYg$ of the Yangian double for any simple finite-dimensional
	Lie algebra ${\bf g}$. Technical computations of the relations in the
	 double are more complicated in general case and we did not make them
		up to the end. Nevertheless in this section we communicate the
 algebraic structure of $\DYg$.

Let $\gg$ be a simple Lie algebra with a standard Cartan matrix
$A=(a_{ij})_{i,j=0}^{r}$, a system of simple roots
$\Pi:=\{\a_1,\ldots,\a_2\}$ and a system of positive roots $\D_+(\gg)$.
Let $e_{i}:=e_{\a_i}$, $h_{i}:=h_{\a_i}$, $f_{i}:=f_{\a_i}:=e_{-\a_i}$,
$(i=1,\ldots,r)$, be Chevalley generators and
$\{e_{\g},f_{\g}\}$, $(\g\in\D)$, be
 Cartan-Weyl basis in $\gg$, normalized so that $(e_{\a},f_{\a})=1$.

Let us first describe $\DYg$ as an algebra.\\
Central extension $\DYg$ of the Yangian $Y(\gg)$ is a formal completion
 (see (\ref{2.7})) of the algebra generated by
the elements $e_{ik}:=e_{\a_i,k}$,
$h_{ik}:=h_{\a_i,k}$, $f_{ik}:=f_{\a_i,k}$,
$(i=1,\ldots,r;\; k=0,1,2,\ldots)$, $c$ and $d$,
 $\mbox{deg}\,e_{ik} = \mbox{deg}\,f_{ik} = \mbox{deg}\,h_{ik} =k$,
 $\mbox{deg}\,c = 0, \mbox{deg}\,d = -1$.
 The relations are written in terms of generating functions
 $$ e^{\pm}_i(u):=\pm\sum_{k\geq
0\atop k<0}e_{ik}u^{-k-1} , \hsp f^{\pm}_i(u):=\pm\sum_{k\geq
0\atop k<0}f_{ik}u^{-k-1}, \hsp h^{\pm}_i(u):=1\pm\h\sum_{k\geq
0\atop k<0}h_{ik}u^{-k-1},$$
 $$e_i(u)=e^+_i(u)-e^-_i(u),\hsp f_i(u)=f^+_i(u)-f^-_i(u):$$
$$[d,e_i(u)]=\frac{\mbox{d}}{\mbox{d}u}e_i(u), \hsp
 [d,f_i(u)]=\frac{\mbox{d}}{\mbox{d}u}f_i(u), \hsp
 [d,h_i^\pm(u)]=\frac{\mbox{d}}{\mbox{d}u}h_i^\pm(u),$$
\begin{eqnarray}
e_i(u)e_j(v)&=&{u-v+\h b_{ij}\over u-v-\h b_{ij}}\ e_i(v)e_j(u) \nn\\
f_i(u)f_j(v)&=&{u-v-\h \bij\over u-v+\h\bij}\ f_i(v)f_j(u) \nn\\
h_i^\pm(u)e_j(v)&=&{u-v+\h\bij\over u-v-\h\bij}\ e_j(v)h_i^\pm(u) \nn\\
h_i^+(u)f_j(v)&=&{u-v-\h\bij-\h\bij c\over u-v+\h\bij-\h\bij c}\
f_j(v)h_i^+(u) \nn\\
h_i^-(u)f_j(v)&=&{u-v-\h\bij\over u-v+\h\bij}\ f_j(v)h_i^-(u) \nn\\
h_i^+(u)h_j^-(v)&=&{u-v+\h\bij\over u-v-\h\bij}\cdot{u-v-\h\bij-
\h\bij c\over u-v+\h\bij-\h\bij c}\
h_j^-(v)h_i^+(u)
\nn\\
{[}e_i(u),f_j(v){]}&=&
{\delta_{i,j}\over \h}\left(\delta(u-(v+\h b_{ii} c))h_i^+(u)-\delta(u-v)
h_i^-(v) \right),
\nn
\end{eqnarray}
\bn
\left\{\begin{array}{l}
{\rm Sym}_{\{k\}}[e_{i}(u_{k_1})[
e_{i}(u_{k_2}) \ldots
[e_{i}(u_{k_{n_{ij}}},
e_{j}(v)]\ldots ]]=0
\\
{\rm Sym}_{\{k\}}[f_{i}(u_{k_1})[
f_{i}(u_{k_2}) \ldots
[f_{i}(u_{k_{n_{ij}}},
f_{j}(v)]\ldots ]]=0
\end{array}\right. \hsp {\rm for}\;\; i\neq j\ ,
\label{Y70}
\ed
Here $\bij = {1\over 2}(\a_i,\a_j)$, $n_{ij} = 1-A_{ij}$.

Hopf algebra $\DYg$ is the double of a Hopf algebra $\widehat{Y}^+(\gg)$
which is isomorphic to  tensor product $Y(\gg)\otimes {\bf C}[c]$
 as an algebra.
The comultiplication in $\Ypg$ is given by the following rules,
where we identify
 elements of $U(\gg)$ with corresponding elements of $Y(\gg)$ generated by
 $e_{i0}, h_{i0}, f_{i0}$.
$$\D(c)=c\otimes 1 +1\otimes c,$$
 $$
\D(x)=x\otimes 1+1\otimes x\ ,\hsp x\in \gg\ ,
$$ $$
\D(e_{i1})=e_{i1}\otimes 1 + 1\otimes e_{i1} +
\h(h_{i0}-c)\otimes e_{i0}-
\h\sum_{\g\in\D_+(g)}f_\g\otimes [e_{\a_i},e_\g]\ ,
$$ $$
\D(f_{i1})=f_{i1}\otimes 1 + 1\otimes f_{i1} + \h f_{i0}\otimes
(h_{i0}+c) +
\h\sum_{\g\in\D_+(g)}[f_{\a_i},f_\g]\otimes e_\g\ ,
$$
\bn
\D(h_{i1})=h_{i1}\otimes 1 + 1\otimes h_{i1} +
\h(h_{i0}-c)\otimes h_{i0}-
\h\sum_{\g\in\D_+(g)}(\a_i,\g)f_\g\otimes e_\g\ .
\label{Y11}
\ed
There is also partial information about comultiplication for basic
 fields in $\DYg$:
$$
\D(e^{\pm}_i(u))=A^\pm\left(e^{\pm}_i(u)\otimes 1
+h_i^{\pm}(u)\otimes e^\pm_i(u)\right)
\quad {\rm mod}\;
	F^\pm\Ypmg\otimes (E^\pm)^2,$$
$$\D(f^{\pm}_i(u))=B^\pm\left(1\otimes f^{\pm}_i(u)
+f^{\pm}_i(u)\otimes h^\pm_i(u)\right)
\quad {\rm mod}\;
(F^\pm)^2\otimes E^\pm\Ypmg,
$$
\bn
\D(h^{\pm}_i(u))= A^\pm
h^\pm_i(u)\otimes h^\pm_i(u)
\quad {\rm mod}\;
F^\pm\Ypmg\otimes E^\pm\Ypmg
\ed
where
$$ A^+ =Id\otimes T_{\h c\otimes 1}^{-1},\hsp
  B^+ =T_{1\otimes \h c}\otimes Id,\hsp A^-=B^-=Id\otimes Id,$$
	 $E^\pm$ and $F^\pm$ are subalgebras generated by the components of
	 $e_i^\pm(u)$ and $f_i^\pm(u)$, $i = 1,2,\ldots ,r=\mbox{rank}\gg$,
		$\Ymg =(\Ypg)^0$.

 Note that, to the contrary to $sl_2$ case, there is no explicit
        formula for the comultiplication in $\Ymg$ except that it is dual
        to the multiplication in $\Ypmg$ with respect to their Hopf pairing.
        This pairing preserves the same decompositions as in $sl_2$ case
         and for the basic fields we have \cite{KT}
 $$<e_i^+(u),f_j^-(v)>= \frac{\delta_{i,j}}{\h (u-v)},\hsp
 <f^+(u),e^-(v)>= \frac{\delta_{i,j}}{\h (u-v)},$$
 $$<h^+(u),h^-(v)>=\frac{u-v+\h\bij}{u-v-\h\bij}\quad\mbox{and}\quad
  <c,d>={1\over \h}.$$
  Analogously to $sl_2$ case, the universal $R$-matrix $\R$ for $\DYg$
  could be reconstructed from that of $DY(\gg)$ by the same procedure.
  Unfortunately,
  complete exact formula for $\R$ is not known in general case.

  {\bf Acknowledgments} The author thanks D.Lebedev and S.Pakuliak
         for fruitful discussions which led to clarifying of the subject
          and to continuation of the work. He thanks for hospitality
          Universite de Reims where this work was started. The work was
          supported by ISF grant MBI300 and Russian Foundation for
                 Fundamental Researches.

\end{document}